# ON THE GEOMETRIC ASPECTS OF ARCHITECTURAL COMPOSITIONS

**Aren Nahapetyan[1], Linda Khachatryan[2]**

[1] *Sapienza University of Rome, Italy, ar.nahapetyan@gmail.com*

[2] *Institute of Mathematics NAS RA, Yerevan, Armenia, linda@instmath.sci.am*

**Abstract:** Different aspects of application of non-Euclidian geometries in architectural compositions are considered. These aspects are illustrated by examples of Armenian (medieval and contemporary) architectural compositions.

> Geometry is the grammar of architecture, and the world around us is the world of geometry.
> Le Corbusier

**Introduction**

It is well known that geometry (of one kind or another) directly influences architectural creativity, opening up vast possibilities for the creation of new forms and compositions.

This paper presents geometric essays outlining elements of well-known geometric theories (Lobachevski's non-Euclidean geometry and Mandelbrot's fractal geometry). We consider different aspects of application of non-Euclidian geometries in architecture. These aspects are illustrated by examples of Armenian (medieval and contemporary) architectural compositions. The paper compiles the series of the previous works [1-3].

In the first part of the work, we consider elements of Lobachevski's non-Euclidean geometry and two of its interpretations (models), the Poincaré model in a half-plane and the Poincaré model in a disc.

Second part is dedicated to fractal geometry. When conducting fractal analysis, various computing tools are used. In this work, the analysis is carried out on the basis of the FrakOut! package, which is very convenient for calculating the fractal parameters of buildings. When evaluating the architectural compatibility of the plan and facade of the temples under consideration, the STATISTICA software package was used to find statistical estimates based on the available data.

**1. Euclidean interpretations of non-Euclidean geometry. Poincaré model**

Contemporary architecture, especially buildings that utilize unconventional design approaches, has always been attractive. Such structures often employ spatial forms inherent to non-Euclidean geometry. From this perspective, certain knowledge of this area of mathematics will be very useful to architects.

In this section, we will provide some elements of non-Euclidean geometry, including the history of its development and solutions to common problems. As for the use of non-Euclidean geometry in Armenian architecture, there are few of them, but they exist.

*Fifth postulate (parallel axiom) of the Euclidean geometry*

Among the axioms of Euclid, the fifth postulate stands alone. While all other postulates are intuitively obvious, their formulations are simple, concise and understandable, the fifth postulate does not have these advantages. Its phrasing is cumbersome, and it takes an effort to understand its meaning: "If a line segment intersects two straight lines forming two interior angles on the same side that are less than two right angles, then the two lines, if extended indefinitely, meet on that side on which the angles sum to less than two right angles."

But the difference between the fifth postulate and others is not limited to this. The fact is that the requirement contained in it refers to the property of a straight line along its entire length. The fifth postulate states that two straight lines are parallel if they do not have a common point, no matter how far they continue. It follows from this that, unlike the other four postulates of Euclid, the fifth postulate cannot be verified in practice.

Naturally, many geometers of those times had a desire to prove the fifth postulate based on the other four. Among them, we note such well-known thinkers as Ptolemy (author of the geocentric system of the world) and Proclus (who wrote the treatise "Comments on Book I of Euclid's Elements"). For example, Proclus found a much simpler equivalent formulation of the fifth postulate (Proclus's axiom): only one straight line parallel to the given one passes through a point outside a straight line.

Apparently, Euclid himself understood that the fifth postulate is special because the proofs of the first 28 theorems of his first book are based only on the first four axioms. The set of those geometric results that follow from the first four axioms of Euclid is called absolute geometry.

In modern times, many attempts to prove the fifth postulate were carried out by the method of contradiction. It is curious to note that in 1763, G.S. Kluge collected and studied about 30 "proofs" of the validity of the fifth postulate and revealed the errors contained in them. From that, he concluded that the fifth postulate is correct and is an independent axiom. Of course, Kluge's conclusion that the fifth postulate is an independent axiom cannot be considered incontestable since, as we know, indirect arguments in favor of any statement cannot be considered as proof of his correctness.

*Non-Euclidian geometry of Lobachevski*

The problem of the fifth postulate of Euclid, which occupied geometers for more than two millennia, was solved practically simultaneously and independently by German mathematician and physicist C.F. Gauss, Hungarian mathematician J. Bolyai and Russian scientist N.I. Lobachevski. Their solution boils down to the fact that the fifth postulate cannot be proved based on other postulates of Euclidean geometry. The new geometry was named Lobachevski's geometry because he was the first who did not hesitate to publicly announce results that were at odds with those that have been accepted for over two thousand years.

Non-Euclidean geometry is a set of statements and theorems that can be obtained on the basis of the system of Euclid's axioms if, instead of the parallel axiom, we accept the opposite statement: on a plane through a point that does not belong to a given straight line, more than one straight line can be drawn that does not intersect with this straight line.

Here are some of the differences between non-Euclidean geometry and Euclidean geometry. In Lobachevski's geometry,

1. There are no similar but unequal triangles. It follows that triangles are equal if their angles are equal. If so, then there is an absolute unit of length, i.e., a segment marked by its properties is similar to a right angle specified by its properties. Such a segment can be, for example, the side of a regular triangle with a given sum of angles.

2. The sum of the angles of any triangle is less than $\pi$ and can be arbitrarily close to zero.

3. Through a point *A*, not lying on a given straight line *a*, there are infinitely many straight lines that do not intersect the straight line *a* and are in the same plane. There are two extreme ones among them, which are called parallel to the straight line *a* in the sense of Lobachevski.

4. If the lines have a common perpendicular, then they diverge infinitely on either side of it.

5. The circumference is not proportional to the radius but grows faster.

It should be noted that the smaller the area in space (on the Lobachevski plane), the less the geometric relations in this area differ from those of the Euclidean geometry. For example, the smaller the triangle, the less the sum of its angles differs from $\pi$; and the smaller the circle, the less the ratio of its length to the radius differs from $2\pi$. A decrease in the area is formally equivalent to an increase in a unit of length. Therefore, with a formal increase in the unit of length, the Lobachevski geometry formulas are transformed into the Euclidean geometry formulas.

The denial of the fifth postulate did not lead to a contradiction with the other four axioms and made it possible to create meaningful and logically harmonious geometry. However, the absence of contradictions does not mean that they cannot appear in the future. Naturally, the problem of consistency of the new, non-Euclidean geometry of Lobachevski arises. To prove the consistency of Lobachevski's geometry, it is necessary to show that it is as consistent as Euclid's geometry, which consistency has not been in doubt for centuries. Such proof was given later when interpretations (models) of non-Euclidean geometry were indicated.

E. Beltrami was the first one who constructed an example of a model called a "pseudo-sphere" (see Fig. 1) in which the first four postulates holding but not the fifth one. From this, it can be seen that non-Euclidean geometry is just as consistent as Euclidean geometry.

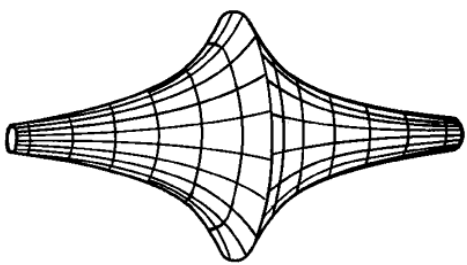

**Fig. 1.** *Pseudo-sphere of Beltrami*

Further, A. Poincaré suggested two models (in a half-plane and in a disc) in which the whole Lobachevski's geometry is presented. Below we describe Poincaré models (we refer to [4-7] for details) and bring some examples of solutions to related problems.

*Poincaré model in a half-plane*

It turns out, that the geometry of Lobachevski is not as clear as the geometry of Euclid. Here, the corresponding models come to the rescue, allowing, based on Euclid's planimetry, to create a meaningful geometric structure in which the axioms of Lobachevski's geometry are fulfilled. Below we present one of them, due to H. Poincaré.

In Lobachevski's planimetry, as in Euclidian, the basic concepts are "point", "straight line", the order relation "lie between" for points of a straight line and "distance between points". All other objects are defined through the above four concepts. The construction of the model begins with the fact that the basic concepts are determined through some notions of Euclidean planimetry, which, in turn, allows defining all other objects in the geometry of Lobachevski. If, in this case, Lobachevski's axioms turn out to be true (verification is carried out using the results of Euclidean geometry), then it is said that the model has been built (implemented). The possibility of constructing such a model testifies the fact that Lobachevski's geometry is as consistent as Euclidian geometry.

The Poincaré model is constructed as follows. A horizontal line is drawn on the Euclidean plane, which divides the plane into two half-planes. This line is called the absolute. Points in the Lobachevski planimetry are the points in the upper half-plane that do not lie on the absolute. Any semicircle centered on the absolute or a ray perpendicular to the absolute, the origin of which lies on the absolute, is called a straight line (see Fig. 2a).

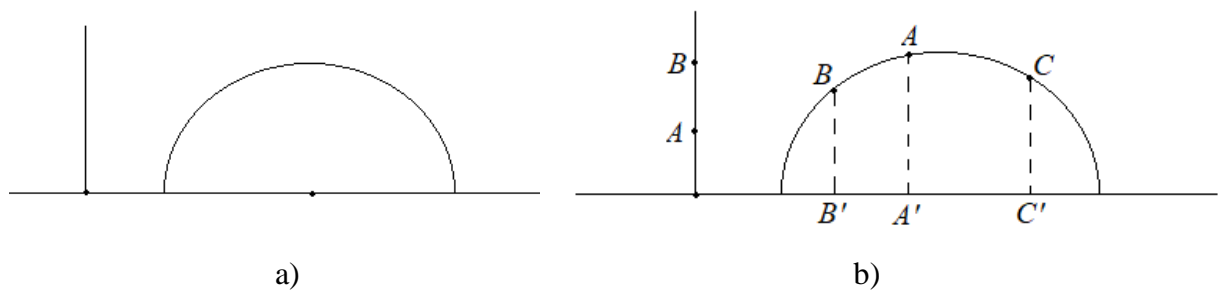

a) b)

**Fig. 2**. *Straight lines and points in the Poincaré model*

Point $A$ lies between two other points $B$ and $C$ on the Lobachevski line if its projection onto the absolute lies between the projections of points $B$ and $C$ (see Fig. 2b). The distance $d(A, B)$ between points $A$ and $B$ is determined by the straight line connecting $A$ with $B$. It can be found in the following way. For the points $A$ and $B$, let the angles α and β be specified as it shown on Fig. 3a. Then

$$d(A,B) = c \left| \ln \frac{|\operatorname{tg}\alpha|}{|\operatorname{tg}\beta|} \right|,$$

where $c$ is some positive constant.

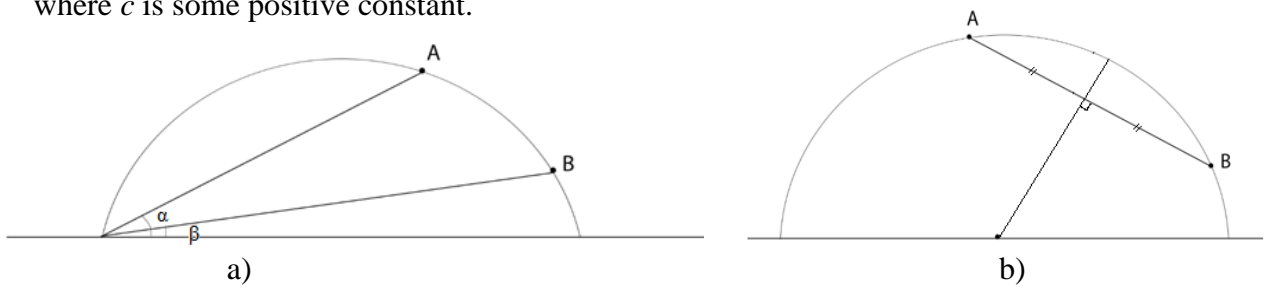

a) b)

**Fig. 3.** *a) Angles α and β in the formula for the distance d(A,B);*
*b) Construction of a semicircle containing two given points A and B*

Let us show that in the constructed model all the axioms of Lobachevski's geometry are fulfilled. To check the first axiom, it is sufficient to construct a semicircle containing the given two points $A$ and $B$ within the framework of Euclidean geometry. This construction is carried out as follows. Let us connect points *A* and *B* with a segment, understood in the sense of Euclidean geometry, and from its center, we lower the perpendicular to the absolute. The point of intersection of this perpendicular with the absolute will be the center of the required semicircle (see Fig. 3b).

The fulfillment of axioms 2, 3 and 4 is obvious. Let us turn to the consideration of the fifth axiom. It is easy to see that an infinite set of straight lines that do not intersect with *a* pass through any point *A* that does not lie on the Lobachevski line *a* (see Fig. 4a). Among these lines, those that have one point in common with *a* on the absolute (lines *b* and *c* on Fig. 4a) stand out. Such Lobachevski lines are called parallel to *a*. We see that there are two straight lines passing through point *A* and parallel to *a*. On Fig. 4b, one can see a pencil of lines parallel to the line *a*, which is perpendicular to the absolute, and passing through a point *A* at infinity lying on the absolute.

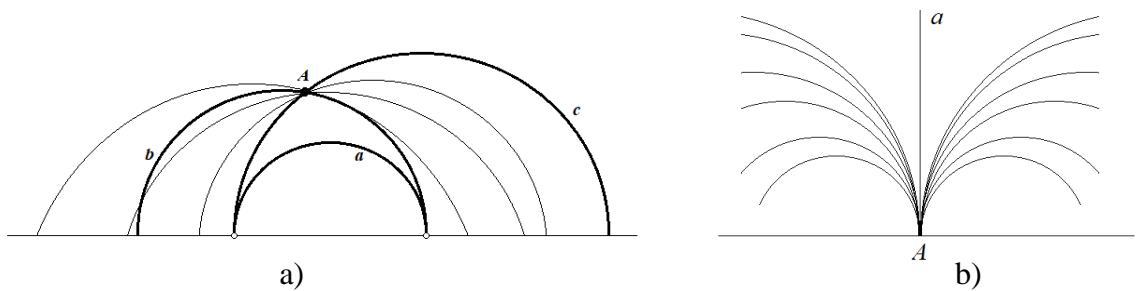

**Fig. 4**. *a) Lines b and c passing through point A and parallel to a;*
*b) Pencil of lines parallel to the line a*

Fig. 5 shows triangles of Lobachevski geometry. It can be seen from this drawing that the sum of their angles is less than π. Note also that the bigger the triangle the more its sum of angles differs from π. To see it, compare triangles *ABC* and *ADE*, as well as triangles *ABC* and *AGF*.

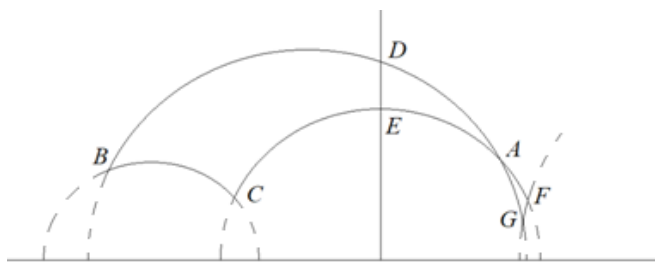 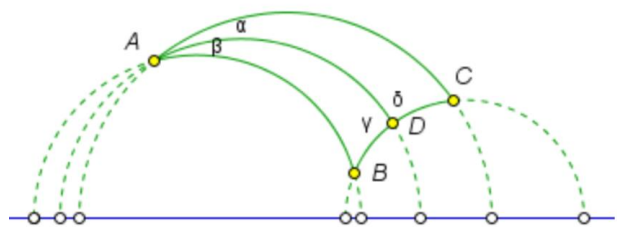

**Fig. 5**. *Triangles of Lobachevski geometry*  **Fig. 6**. *Triangles ABD and ACD*

When solving problems of Lobachevski planimetry (in the Poincaré model), the conditions of the problem are translated into the language of Euclidean geometry, after which the problem is solved by its methods. Let us give some examples.

*Problem.* Show that triangles in the Lobachevski plane can have different sums of angles.

*Solution.* Consider the triangles *ABD* and *ACD* shown on Fig. 6. Let us calculate the sums $\sigma_{ABD}$ and $\sigma_{ACD}$ of their angles. We have
$$\sigma_{ABD} + \sigma_{ACD} = (\beta + \gamma + \angle B) + (\alpha + \delta + \angle C).$$
Using the facts that $\alpha + \beta = \angle A$ and $\gamma + \delta = \pi$, further we can write
$$\sigma_{ABD} + \sigma_{ACD} = \angle A + \angle B + \angle C + \pi = \sigma_{ABC} + \pi.$$
Hence, one of the sums, $\sigma_{ABD}$ or $\sigma_{ACD}$, is less than π. Suppose, $\sigma_{ACD} < \pi$. Then from the equality
$$\sigma_{ABD} + (\sigma_{ACD} - \pi) = \sigma_{ABC},$$
it follows that $\sigma_{ABD} > \sigma_{ABC}$. Thus, we have found two triangles, *ABD* and *ABC* with different sums of its angles.

*Problem.* Let *A*, *B* and *C* be points in the Lobachevski plane, not lying on one straight line, and let the straight line *a* does not pass through any of these points. Prove that if the straight line a passes through the point of the segment *AB*, then it passes through the point of the segment *AC* or *BC*.

*Solution.* In Euclidean geometry, the following statement is true: two circles intersect if and only if one of them passes through the inner point of the other circle. Let the line *a* (in the sense of Lobachevski) intersect the segment *AB* (see Fig. 7). Then either point *A* will be internal for the Euclidean circle *a*, or that will be point *B*. Suppose that this is point *B*. Then if point *C* is an internal point for the Euclidean circle *a*, then the line *a* intersects the segment *AC*. Otherwise, *a* intersects the segment *BC*.

Further, we will prove the analogue of the Pythagorean theorem in the Lobachevski's geometry.

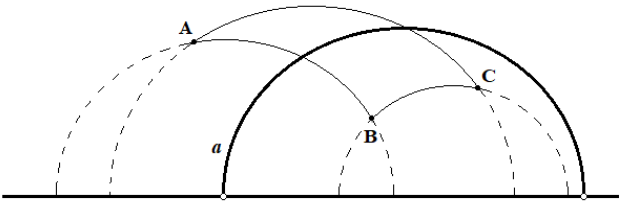
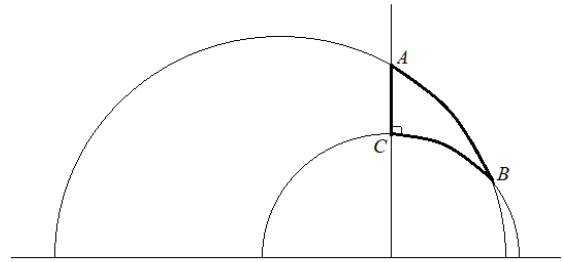

**Fig. 7**. *Line a intersecting segment AB*      **Fig. 8**. *Right triangle ABC*

**Theorem.** *For any right triangle, the following equality holds*
$$\operatorname{ch} c = \operatorname{ch} a \cdot \operatorname{ch} b,$$
*where c represents the length of the hypotenuse and a and b the lengths of the triangle's other two sides.*

**Proof.** Consider the right triangle *ABC* with the right angle *C* (see Fig. 9). Without loss of generality, we will assume that the vertices *A*, *B* and *C* of this right triangle correspond to complex numbers $ri$, $u + vi$ and $I$, respectively, where $r > 1$ and $u^2 + v^2 = 1$, since this can always be achieved with some non-Euclidean motion. Then we have

$$\operatorname{ch} a = \operatorname{ch}|AC| = \operatorname{ch} d(A, C) = 1 + \frac{|ri - i|^2}{2r} = \frac{1 + r^2}{2r},$$

$$\operatorname{ch} b = \operatorname{ch}|BC| = \operatorname{ch} d(B, C) = 1 + \frac{|u + vi - i|^2}{2v} = \frac{1 + u^2 + v^2}{2v} = \frac{1}{v},$$

$$\operatorname{ch} c = \operatorname{ch}|AB| = \operatorname{ch} d(A, B) = 1 + \frac{|ri - u - vi|^2}{2rv} = \frac{u^2 + v^2 + r^2}{2rv} = \frac{1 + r^2}{2rv}.$$

Hence,

$$\operatorname{ch} a \cdot \operatorname{ch} b = \frac{1 + r^2}{2r} \cdot \frac{1}{v} = \frac{1 + r^2}{2rv} = \operatorname{ch} c.$$

□

Lobachevski noticed that the non-Euclidean geometry he created in an infinitesimal, that is, in the first approximation, coincides with the geometry of the Euclidean plane. Let us illustrate this with the example of the Pythagorean theorem. Using the series expansion of the hyperbolic cosine

$$\operatorname{ch} z = 1 + \frac{z^2}{2!} + \frac{z^4}{4!} + \frac{z^6}{6!} + \cdots,$$

we obtain for the Pythagorean theorem the following relation

$$1 + \frac{c^2}{2!} + \frac{c^4}{4!} + \cdots = \left(1 + \frac{a^2}{2!} + \frac{a^4}{4!} + \cdots\right)\left(1 + \frac{b^2}{2!} + \frac{b^4}{4!} + \cdots\right).$$

Keeping the lower order terms now, we arrive at the Pythagorean theorem of Euclidean geometry: $c^2 = a^2 + b^2$.

*Poincaré disc model*

In addition to the above model, Poincaré proposed the model of the Lobachevski plane in a disc. This model is described as follows. The plane itself is represented as the inner part of the circle of radius 1 with the center at the origin of the complex plane, that is, $\mathcal{B} = \{z \in \mathbb{C}: |z| < 1\}$. The points of the circumference bounding the disc do not belong to the plane, and the circumference itself is called the absolute.

This model can be obtained from Poincaré's model in a half-plane through the following transformation

$$w = \frac{z - i}{z + i}.$$

Indeed,

$$|w| = \left|\frac{x + i(y - 1)}{x + i(y + 1)}\right| = \sqrt{\frac{x^2 + (y - 1)^2}{x^2 + (y + 1)^2}} \leq 1.$$

The inverse transition is obviously carried out by means of the mapping

$$w = i\frac{1 + z}{1 - z}.$$

The lines in this model are the arcs of the circles perpendicular to the absolute. Recall that two curves are called perpendicular if their tangents at the point of intersection of the curves are mutually perpendicular. The lines forming a diameter of the absolute are called special.

It is quite obvious (see Fig. 9a) that in the defined model, if a straight line *a* and a point *A* not belonging to it are given, then through the point *A* it is possible to draw an arbitrarily large number of straight lines that do not intersect the straight line *a*. From Fig. 9a it can be seen that among such lines, two stand out – line *b* and line *c*. They are characterized by the fact that they have a common point with line *a* on the absolute. It is these two straight lines in the Poincaré model that are called parallel to the straight line *a*. All others passing through the point *A* and having no common points with it are inside the angle formed by straight lines *b* and *c*. Such straight lines are called diverging.

It follows from what has been said that in Lobachevski's geometry, parallel lines do not possess the property of transitivity, i.e., from the fact that, for example, line *b* is parallel to line *a*, and line *a* is parallel to line *c*, it does not follow that line *b* is parallel to line *c* (see Fig. 9a). However, this fact takes place in Euclidean geometry.

The choice of one or the other model depends on the task. In particular, the fact that the sum of the angles of a triangle in Lobachevski's geometry is less than $\pi$ is much easier to see on a disc model (see Fig. 9b).

Let us show how the distance between points is determined in this Poincaré model of the Lobachevski's plane. Let points *A* and *B* be given, and let the straight line (in the sense of Lobachevski) connecting these points intersects the absolute at points *C* and *D* (see Fig. 9c). Then the distance $d(A, B)$ between points *A* and *B* is determined as follows:

$$d(A, B) = \left| \ln \frac{CB/CA}{BD/AD} \right|,$$

where *CB*, *CA*, *BD* and *BA* are the lengths of the corresponding arcs in the sense of Euclid.

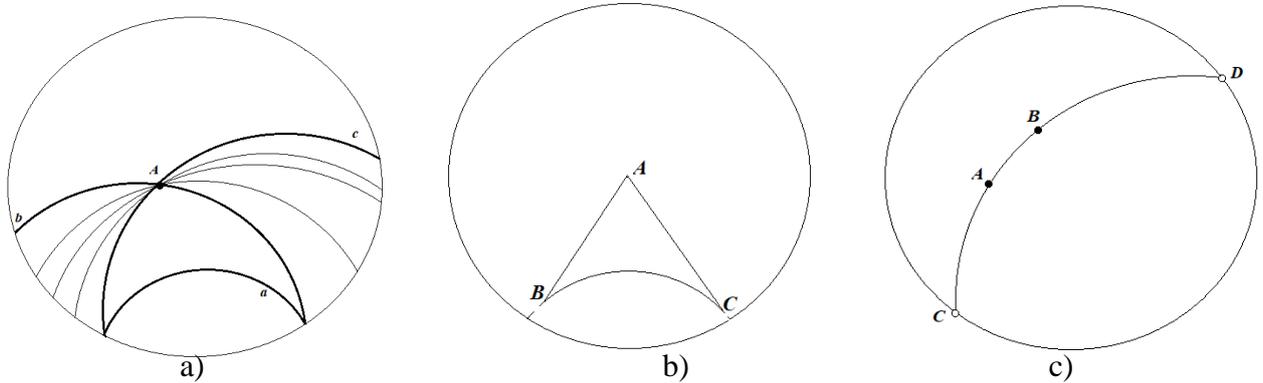

a)  b)  c)

**Fig. 9**. *Poincaré disc model: a) lines drawn through point A and parallel to a;*
*b) triangle; c) distance*

*Problem.* Let points *A* and *B* lie on the same singular line in the Poincaré disc model, and let their Euclidean distances to the center of the disc be known. Calculate the distance from *A* to *B* in the sense of Lobachevski's geometry.

*Solution.* Without loss of generality, let us assume that points *A* and *B* belong to the diameter lying on the *OY* axe of the complex plane.

Let the Euclidean distances to the center from points *A* and *B* be equal to α and β, respectively, $0 \leq \alpha < \beta < 1$. Then to the point *A* it corresponds the complex number $i\alpha$, and to the point *B* it corresponds the number $i\beta$. Thus, we can write

$$d(A,B) = d(i\alpha, i\beta) = 2\mathrm{arth}\left|\frac{i\alpha - i\beta}{1 - i\alpha\overline{(i\beta)}}\right| = 2\mathrm{arth}\left|\frac{\alpha - \beta}{1 - \alpha\beta}\right| = \ln\left(\frac{1 + \left|\frac{\alpha - \beta}{1 - \alpha\beta}\right|}{1 - \left|\frac{\alpha - \beta}{1 - \alpha\beta}\right|}\right) =$$

$$= \ln\frac{(1+\beta)(1-\alpha)}{(1-\beta)(1+\alpha)} = \ln\frac{1-\alpha}{1-\beta} + \ln\frac{1+\beta}{1+\alpha} = \ln\left(1 + \frac{\beta - \alpha}{1-\beta}\right) + \ln\frac{1+\beta}{1+\alpha}.$$

*Elements of Lobachevski geometry in architecture*

Lobachevski's geometry is applied indirectly to architecture through models, such as the pseudo-sphere model or the Poincaré model in a disc. Using these models, architects can design complex non-Euclidean forms, which in modern buildings are represented as concave or saddle-shaped.

Lobachevski's geometry inspired many architectures of a new time as well. Its advances allowed the development of new original design methods for the new architecture. Among the architects of modern times, we will name Ernst Neufert, Zaha Hadid, Frank Gehry, and the great architect Antoni Gaudi. Examples of buildings constructed with the usage of non-Euclidian geometry include Casa Batllo, built in 1900 by famous architect Antonio Gaudi in Art Nouveau style, where main idea is to mimic forms that are found in the nature, which gives building a unique look; Cleveland Clinical Center for Medical Art and Photography, built by Frank Gehry in postmodern style; Shell House, built in 2008 by Japanese architecture studio ARTechnic architects (see Fig. 10).

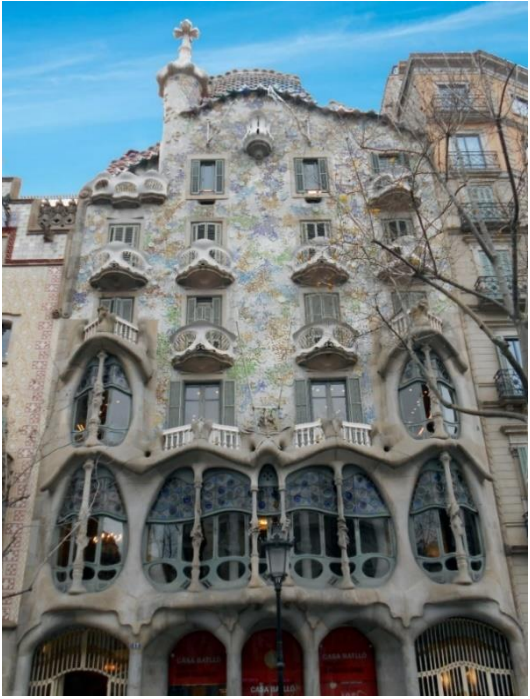

*a) Casa Batllo*

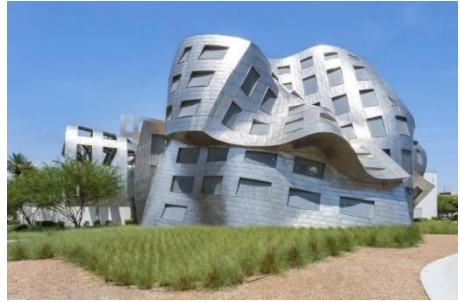

*b) Cleveland Clinical Center for Medical Art and Photography*

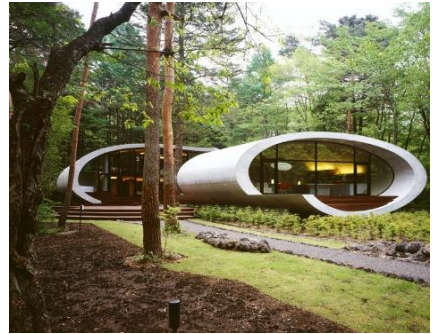

*c) Shell House*

**Fig. 10.** *Architecture inspired by Lobachevski's geometry*

Triangles of Lobachevski geometry, called sails, are also widely used. For example, they are used in architecture for the construction of domes. These triangles are used to solve the problem of how to put a dome with a round base on a square structure.

One can find a lot of examples of mentioned construction in Armenian churches. Fig. 11a shows Haghpat monastery, which is a medieval monastery complex in Haghpat, Armenia, built between the 10th and 13th century. On the right-hand side of Fig. 11a, triangles of Lobachevski geometry are clearly viewed, which are used to withhold the dome. It should be noted that apart from one or two minor restorations carried out in the eleventh and twelfth centuries, the church has retained its original character. The monastery has been damaged many times. Sometime around 1130, an earthquake destroyed parts of Haghpat monastery and it was not restored until fifty years later. It also suffered numerous attacks by armed forces in the many centuries of its existence and from a major earthquake in 1988. Nevertheless, much of the complex is still intact and stands today without substantial alterations.

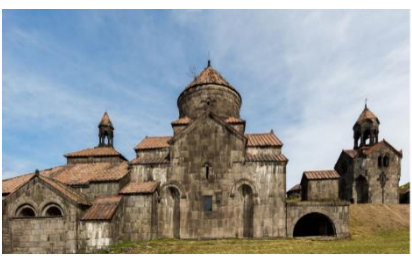 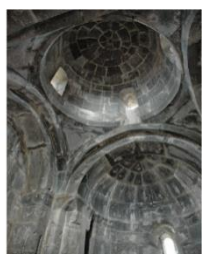 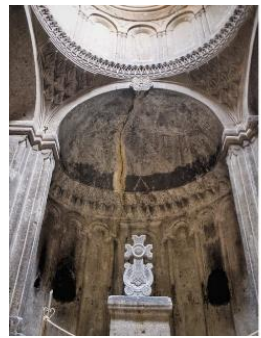 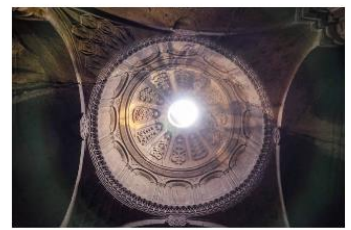

a)          b)

**Fig. 11**. *Sails in Armenian architecture: a) Haghpat monastery; b) Geghard Church*

Another example (see Fig. 11b) is the dome in the Geghard church. This is a medieval monastery in the Kotayk province of Armenia, being partially carved out of the adjacent mountain. Here, the sails are used not only as a part of construction but also as a decoration.

## 2. Fractal geometry and quantitative evaluation of the aesthetic appeal of architecture compositions

The basic figures of classical (Euclidean) geometry are simple and clear: circle, sphere, cylinder, pyramid, etc. Their surfaces assumed to be ideal, both in terms of their shapes and in terms of their surface properties. Given such limitations, classical geometry can only describe a very narrow class of natural structures and phenomena but not complex physical objects, such as the shapes of clouds and mountains, tree crowns, or the human bronchial system.

Geometry describing non-standard forms was proposed by B. Mandelbrot [8] based on the concept of fractal introduced by him. It was called fractal geometry. In contrast to classical Euclidean geometry, where objects are idealized (their surfaces are assumed to be perfectly smooth, without any irregularities, cracks, or breaks), the fractal geometry studies the patterns inherent in natural objects, processes, and phenomena with the presence of roughness, brokenness, and other complexities (see, for example, [9-11]). It offers a variety of ways to describe and measure both natural and man-made objects.

*Fractals and fractal geometry*

According to Mandelbrot [8], a fractal is a structure consisting of parts that are in some sense similar to the whole (or to each other). From a mathematical point of view, a fractal is a geometric figure (a set of points in Euclidean space) whose fractal dimension (the Hausdorff-Besicovitch dimension) is either fractional or exceeds its topological dimension.

The fractal dimension for a finite set $G$ in $\mathbb{R}^n$ is defined as follows. Let $\mathbb{Z}^n(\Delta)$ be a cubic lattice in $\mathbb{R}^n$ with cube (cell) edge length equal to $\Delta$. Let $N(\Delta)$ be the minimum number of lattice cubes needed to cover $G$. Then the fractal dimension $D$ of $G$ is defined by

$$D = - \lim_{\Delta \to 0} \frac{\ln N(\Delta)}{\ln \Delta}.$$

This dimension can be defined equivalently based on the following requirement:

$$\lim_{\Delta \to 0} N(\Delta)\Delta^d = \begin{cases} 0 & \text{if } d > D, \\ \infty & \text{if } d < D. \end{cases}$$

Thus, the dimension $D$ of the set $G$ is essentially the boundary that shows that if $d < D$, then the number of cubes $N(\Delta)$ is insufficient to cover the set $G$, and if $d > D$, then the number of cubes is excessive for coverage.

It is generally accepted that the fractal dimension is a characteristic property of fractals, i.e., if the dimension $D$ is not an integer, then the set $G$ is considered a fractal. The Hausdorff-Besicovitch dimension increases with the degree of tortuosity of the object. For a straight line, it is equal to one, for a slightly tortuous line it is 1.03, for a more tortuous line it is 1.16, and for a strongly tortuous line it is 1.57, and so on.

In addition to natural fractals (such as island coastlines, snowflakes, crystals, and heads of broccoli), there are also artificial (non-natural) fractals. The first examples of non-natural fractals were constructed at the end of the nineteenth century in connection with purely mathematical problems of function theory. From the point of view of classical mathematical analysis, they had extremely unusual properties. For example, this is the Cantor set (Cantor dust), the nowhere differentiable Weierstrass function, the Koch snowflake, the Brownian curve on the plane, etc. For some of them, fractal dimensions have been calculated: the Cantor set has a fractal dimension $D = \ln 2 / \ln 3$, and for a Brownian curve on a plane, it is equal to 2, that is, exceeds its topological dimension.

Practical methods for calculating the fractal dimension play an important role. One of the most popular methods is the method of counting cells that have a non-empty intersection with the image being studied (box-counting dimension method). Apparently, W. Lorenz [12] and C. Bovill [13] were the first to study and use this method most fully. Let us describe in general terms the algorithm for applying this method.

In the first step, a cubic (square) grid with the cell edge length (scale) equal to $\Delta$ is superimposed on the image under study. Initially, $\Delta$ is taken to be equal to $L$, where $L$ is the length of the rectangle containing all the images. Let $N(\Delta)$ be the number of all cubes that have a non-empty intersection with the image under study. Next, the ration $-\log N(\Delta) / \log \Delta$ is considered, and its behavior is investigated under stepwise changes in the scale $\Delta$.

The scale is reduced by half at each step. The process can continue indefinitely, but in practical applications, it is stopped depending on the requirements of the task. The slope of the graph of $\log N(\Delta)$ from $-\log \Delta$ gives an approximate value of the fractal dimensions of the image.

*Fractals in architecture*

In architecture, fractal principles are used in the design of objects using computer modeling. These principles can be used to create unique and very interesting architectural forms (see, for example, [11]).

Another application of fractal geometry in architecture came from the psychology where objective (quantitative) methods of the aesthetic appeal evaluation of architectural compositions where developed. The first studies in the psychology of fractal perception were conducted by J. Sprott's group [14]. Similar studies were conducted by R. Taylor [15], who demonstrated that the overwhelming majority of experimental participants preferred fractal patterns. Sprott found that the attractiveness of fractal objects correlates with their fractal dimension. The results of experiments showed that participants preferred objects with a fractal dimension between 1.1 and 1.5.

Research on the psychology of fractal perception is currently underway (see, for example, [16-20]). These studies have confirmed the hypothesis that the aesthetic appeal of architectural compositions is largely determined by the value of their fractal dimension.

There are many monuments of world architecture for which fractal analysis has been carried out. These are remarkable Gothic cathedrals in Europe, beautiful mosques of Islamic architecture, unique Hindu temples (see, for example, [21-24]).

Regarding Armenian architecture, publications on fractal analysis of its compositions are nearly absent. In a recently published works [2-3], the author conducted a fractal analysis of three major monuments of medieval Armenia (the Hripsime temple, the Zvartnots Cathedral and the Cathedral of the Holy Virgin in Ani) as well as some modern compositions (the Yerevan Cascade, the Government House #2 on the Republic Square, the Cathedral of Saint Gregory the Illuminator, and the Church of the Holy Trinity). In this section we briefly present the obtained results.

*Fractal analysis of medieval Armenian temples*

In [2], a fractal analysis of the temples of Hripsime and Zvartnots, as well as the Ani Cathedral was presented. A statistical analysis of the compatibility of the plan and facade of these buildings was also carried out.

The obtained results are summarized in Table 2. They show that the temples under consideration have high architectural attractiveness, and their plan and facade are in excellent agreement with each other.

**Table 2.** Calculations of the fractal dimension of the facade and plan of Armenian temples

| Calculation of fractal dimension between: | | fractal dimension | | | | | |
|---|---|---|---|---|---|---|---|
| | | Hripsime | | Zvartnoc | | Ani cathedral | |
| large grid size | small grid size | facade | plan | facade | plan | facade | plan |
| 200 | 100 | 1.46 | 1.74 | 1.64 | 1.67 | 1.56 | 1.48 |
| 100 | 50 | 1.48 | 1.58 | 1.54 | 1.57 | 1.53 | 1.50 |
| 50 | 25 | 1.49 | 1.49 | 1.48 | 1.49 | 1.56 | 1.43 |
| 25 | 12.5 | 1.49 | 1.51 | 1.47 | 1.43 | 1.5 | 1.13 |
| average fractal dimension | | 1.48 | 1.58 | 1.533 | 1.540 | 1.537 | 1.385 |

**Hripsime.** The temple was built by Catholicos Komitas in 618 to the east of Echmiadzin on the burial site of Saint Hripsime. It is a central-domed structure with an internal cross-shaped base. It is a recognized masterpiece of Armenian architecture.

The temple of Hripsime has an average fractal dimension of 1.48. The calculations also show that the standard deviation of these data from the average is 0.014. Regarding the architectural plan, the following estimates were obtained: the average fractal dimension is 1.58 with a standard deviation of 0.113. The correlation between the fractal dimensions of the facade and the plan is – 0.997.

**Zvartnots.** The Zvartnots Cathedral was founded by Catholicos Nerses III in the middle of the 7th century, not far from Vagharshapat (Echmiadzin) in the place where, according to legend, Gregory the Illuminator and the king of Armenia Trdat met. This majestic temple is a tetraconch (a central-domed structure with a plan in the form of a cross with rounded ends).

The Zvartnots temple has an average fractal dimension of 1.533 with standard deviation 0.008 of obtained data from the average. For the architectural plan, the average fractal dimension is 1.54 with a standard deviation of 0.104. The correlation between the fractal dimensions of the facade and the plan is 0.974.

**Cathedral of the Holy Virgin in Ani.** The Ani Cathedral is the pinnacle of Armenian architecture of the 9th-11th centuries. It is a prototype of Gothic architecture. Its architectural forms are similar to European Gothic.

Regarding Gothic, we note that there is a very reasonable assumption that the first object where Gothic principles were applied was not the Cathedral of Saint-Denis (a suburb of Paris), but the Cathedral of the Holy Virgin in Ani. The interior of this temple clearly contains such architectural compositions as elongated pointed arches, bunches of columns with ribbed vaults. These compositions were developed in Gothic architecture, which was widespread in Western Europe.

In his major work [25], Professor of the University of Vienna J. Strzygowski writes: "Consequently, it remains to be recognized that the Armenians built in the Gothic style approximately 150 years earlier than was the case in Europe".

The Ani Cathedral has an average fractal dimension of 1.537 with standard deviation 0.029. The architectural plan has an average fractal dimension 1.385 with a standard deviation of 0.172. The correlation between the fractal dimensions of the facade and the plan is 0.797.

*Fractal analysis of contemporary Armenian architectural structures*

In [3], a fractal analysis was conducted for several contemporary Armenian architectural structures with obvious fractal motifs: the Yerevan Cascade, the Government House #2 on the Republic Square, the Cathedral of Saint Gregory the Illuminator, and the Church of the Holy Trinity. The obtained results are summarized in Table 3. Below we will present the fragments of the process of calculating the fractal dimension for mentioned constructions.

**Table 3.** Calculations of the fractal dimension of modern Armenian architectural structures

| Calculation of fractal dimension between: | | Cascade | Government House #2 | Cathedral of St. Gregory the Illuminator | Holy Trinity Church |
|---|---|---|---|---|---|
| large grid size | small grid size | | | | |
| 128 | 64 | 1.51 | 1.56 | 1.66 | 1.63 |
| 64 | 32 | 1.49 | 1.51 | 1.60 | 1.63 |
| 32 | 16 | 1.44 | 1.52 | 1.52 | 1.59 |
| 16 | 8 | 1.38 | 1.49 | 1.46 | 1.52 |
| average fractal dimension | | 1.455 | 1.52 | 1.56 | 1.593 |

**Yerevan Cascade**. The Yerevan Cascade is an architectural and monumental complex consisting of a multi-level staircase with terraces, fountains, sculptures, and exhibition halls. Its construction began in 1980, but was halted in the late 1980s due to the Spitak earthquake, the collapse of the USSR, and the First Karabakh War. The architecture combines elements of Soviet modernism with traditional motifs of Armenian stone architecture.

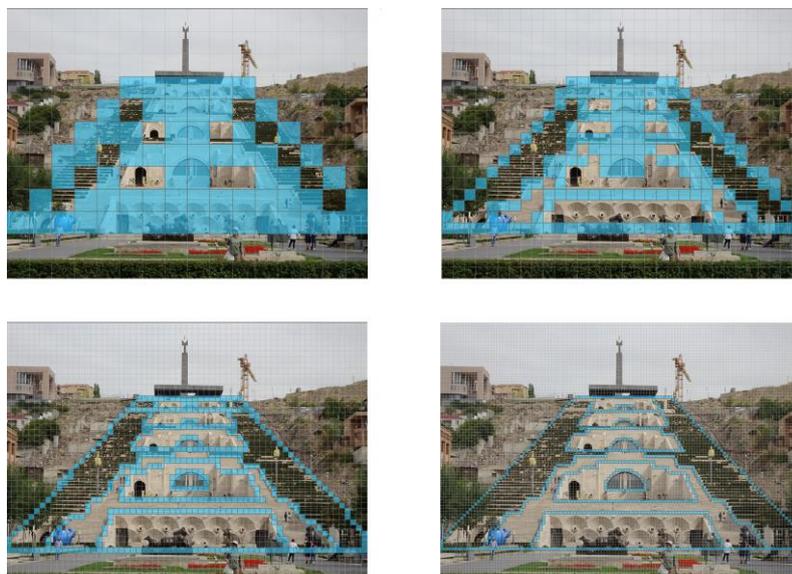

**Fig. 12.** *Calculation of the fractal dimension of the facade of the Yerevan Cascade*

Fig. 12 shows fragments of the process of calculating the fractal dimension of the Cascade. From the obtained data it follows that the Cascade has an average fractal dimension of 1.455. The calculations also show that the standard deviation of these data from the average is 0.058.

**Government House #2 on the Republic Square.** Republic Square is the central square of Yerevan. Government House #2 on this square was built in 1955 according to the design of Samvel Safaryan, Rafael Israelyan, and Varazdat Arevshatyan. From 1996 to 2016, the building belonged to the Ministry of Foreign Affairs of the Republic of Armenia. Fig. 13 shows fragments of the process of calculating the fractal dimension for this building. Tthe Government House has an average fractal dimension of 1.52 with standard deviation 0.029 of these data from the average.

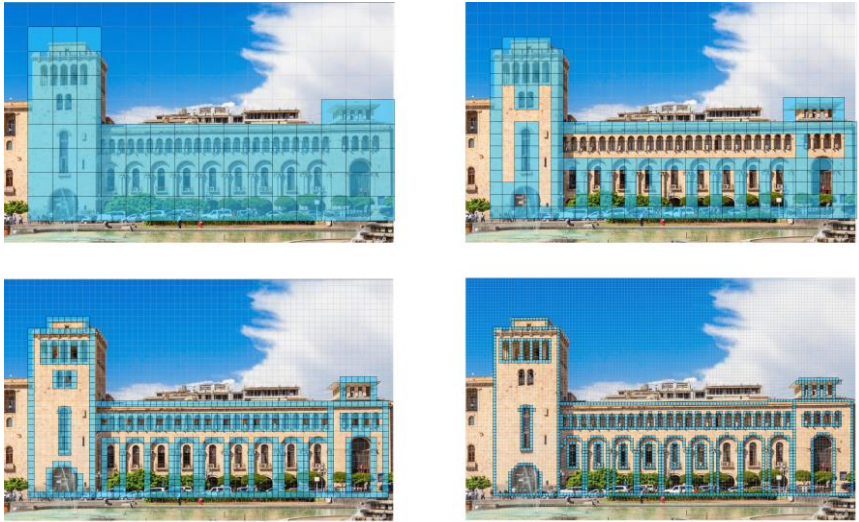

**Fig. 13.** *Calculation of the fractal dimension of the facade of the Government House #2*

**Cathedral of St. Gregory the Illuminator.** The Cathedral of Saint Gregory the Illuminator is the largest cathedral built in Yerevan. It was erected to commemorate the 1700th anniversary of Armenia's adoption of Christianity as the state religion. Construction began in 1997 and took approximately four years. The cathedral's architect is Stepan Kyurkchyan.

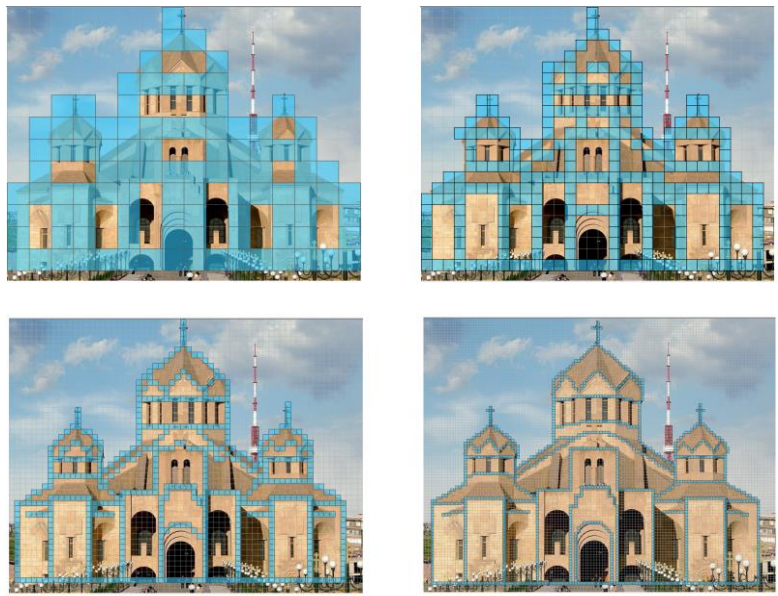

**Fig. 14.** *Calculation of the fractal dimension of the facade of the St. Gregory the Illuminator cathedral*

The Cathedral of St. Gregory the Illuminator has an average fractal dimension of 1.56 with standard deviation 0.088. Fig. 14 shows fragments of the process of calculating its fractal dimension.

**Church of the Holy Trinity.** The Holy Trinity Church was built in 2003 in Yerevan, in the Malatia--Sebastia district. The design was developed by honored architect Baghdasar Arzumanyan. The church is modeled after the Zvartnots Cathedral.

Fig. 15 shows fragments of the process of calculating the fractal dimension of the church. The results show that the Church of the Holy Trinity has an average fractal dimension of 1.593 with standard deviation 0.052.

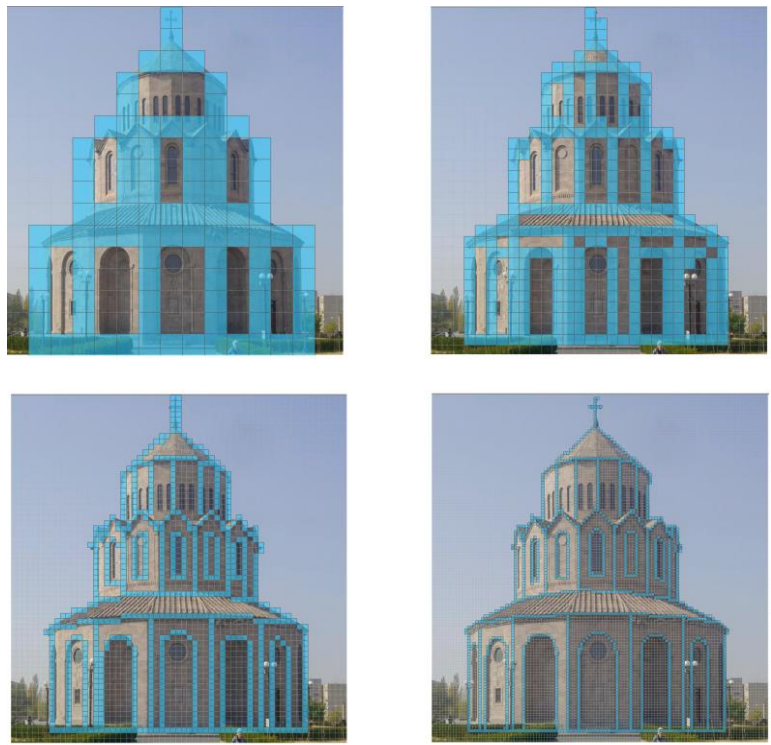

**Fig. 15.** *Calculation of the fractal dimension of the facade of the Holy Trinity Church*

**Conclusion**

Elements of non-Euclidean geometry are often found in Armenian architectural compositions. Fractal analysis of the examined Armenian medieval temples as well as contemporary architecture showed a high level of consistency between subjective and objective assessments of their aesthetic appeal.

**References**

[1]. A. Nahapetyan, Euclidean interpretations of non-Euclidean geometry. Poincaré model. "Quantum" college, Yerevan, 2022 (unpublished essay).
[2]. A. Nahapetyan, Fractal geometry and quantitative evaluation of the aesthetic appeal of ancient Armenian architecture monuments. Mathematical Problems of Computer Science, 63, 2025, 42-53. Doi: https://doi.org/10.51408/1963-0130
[3]. A. Nahapetyan, Elements of fractal geometry in Armenian architectural compositions. Izvestia NAS RA, Mathematics, 2025 (accepted).
[4]. R. Courant, H. Robbins, What is mathematics? Oxford University Press, 2nd edition, 1996.

[5]. B.F. Kagan, The foundations of geometry (in Russian). State Publishing House of Technical and Theoretical Literature, 1949.

[6]. D. Marshall, P. Scott, A brief history of non-Euclidean geometry. Australian Mathematics Teacher, 60 (3), 2004, 2-4.

[7]. V.V. Prasolov, Geometry of Lobachevski (in Russian). MCCMI, Moscow, 3rd ed., 2004

[8]. B.B. Mandelbrot, The Fractal Geometry of Nature, Freeman, San Francisco, 1982.

[9]. R.M. Crownover, Introduction to Fractals and Chaos, Jones and Bartlett Publishers, 1995.

[10]. J. Feder, Fractals, Publisher New York: Plenum Press, 1988.

[11]. N.A. Salingaros, Architecture, Patterns, and Mathematics. Nexus Network Journal, 1 (1), 1999, 75-86. Doi: https://doi.org/10.1007/s00004-998-0006-0

[12]. W. E. Lorenz, Fractals and Fractal Architecture: PhD thesis, Vienna University of Technology, Vienn, 2003.

[13]. C. Bovill, Fractal Geometry in Architecture and Design, Birkhäuser, 1996, 73-92.

[14]. J.C. Sprott, Strange attractors: Creating patterns in chaos. New York, M&T Books, 1993.

[15]. R.P. Taylor, Reduction of physiological stress using fractal art and architecture. Leonardo, 39 (3), 2006.

[16]. S. Pyankova, Fractal analysis in psychology: perception of self-similar objects (in Russian). Psychological Studies, 9 (46), 2016, 14 pp.

[17]. S. Pyankova, Subjective estimates of visual complexity and aesthetical appeal of fractal images: individual differences and genetic influences (in Russian). Psychological Studies, 12 (63), 2019, 16 pp.

[18]. N.A. Salingaros, Fractal art and architecture reduce physiological stress. Computers and Graphics, 27 (5), 2003, 813-820.

[19]. C.M. Hagerhall, T. Purcell, R. Taylor, Fractal dimension of landscape silhouette outlines as a predictor of landscape preference. Journal of Environmental Psychology, 24 (2), 2004, 247-255.

[20]. F.I. Mavrikidi, Fractal mathematics and the nature of change (in Russian). Delphis, 54 (2), 2008.

[21]. P.K. Acharya, An encyclopaedia of Hindu architecture. Oxford University Press, 2010.

[22]. H.D.A. Ismail, Fractal analysis of masterpieces of Islamic architecture – Ahmad Shah Mosque and Taj Mahal: justification of the method and experience of application (in Russian). Architecture and Modern Information Technologies, 4 (21), 2012, 11 pp.

[23]. H.D.A. Ismail, M.Yu. Shishin, Application of multi-stage fractal method in the analysis of the masterpiece of Islamic architecture – Ahmad Shah Mosque (in Ruassin). Art of Eurasia, 3 (10), 2018, 37-47. Doi: https://doi.org/10.25712/ASTU.2518-7767.2018.03.003

[24]. I.A. Mayatskaya, G.I. Fazylzyanova, B M. Yazyev, S.B. Yazyeva, Fractality of gothic architecture (in Russian). Engineering Bulletin of the Don, 4, 2025, 15 pp.

[25]. J. Strzygowski, Die Baukunst der Armenier und Europa: Ergebnisse einer vom Kunsthistorischen Institute der Universitat Wien 1913 durchgefuhrten Forschungsreise. Kunstverlag Anton Schroll & Co., 1918.